\documentclass[10pt,a4paper,oneside]{amsart}

\usepackage{amsfonts}

\usepackage{amssymb}

\usepackage{latexsym}

\usepackage{eufrak}

\usepackage{euscript}

\usepackage{amscd}

\usepackage{url}

\mathchardef \poltheta = "0112

\allowdisplaybreaks

\newtheorem{theorem}{Theorem}[subsection]
\newtheorem{conjecture}[theorem]{Conjecture}
\newtheorem{proposition}[theorem]{Proposition}
\newtheorem{lemma}[theorem]{Lemma}
\newtheorem{corollary}[theorem]{Corollary}

\newtheorem*{theorem*}{Theorem}
\newtheorem*{proposition*}{Proposition}
\newtheorem*{lemma*}{Lemma}
\newtheorem*{corollary*}{Corollary}
\newtheorem*{definition*}{Definition}

\theoremstyle{remark}

\newtheorem*{remark*}{Remark}
\newtheorem*{example*}{Example}

\newcommand{\CC}{\mathbb{C}}

\newcommand{\FF}{\mathbb{F}}

\newcommand{\HH}{\mathbb{H}}
\newcommand{\NN}{\mathbb{N}}

\newcommand{\RR}{\mathbb{R}}

\newcommand{\ZZ}{\mathbb{Z}}

\newcommand{\calX}{\mathcal{X}}

\newcommand{\Id}{\mathbf{1}}
\newcommand{\idest}{\textit{i.e.} }

\newcommand{\Rh}[2]{[#1 \ #2]}

\newcommand{\set}[2]{\left\{#1\,\mid\,#2\right\}}
\newcommand{\tp}[1]{{^{t}}\!#1}

\DeclareMathOperator{\Char}{char}
\DeclareMathOperator{\Cof}{Cof}

\DeclareMathOperator{\Disc}{Disc}
\DeclareMathOperator{\Gal}{Gal}
\DeclareMathOperator{\GL}{\mathrm{GL}}

\DeclareMathOperator{\Jac}{Jac}

\DeclareMathOperator{\PGL}{\mathrm{PGL}}
\DeclareMathOperator{\Res}{Res}
\DeclareMathOperator{\Trans}{Trans}
\DeclareMathOperator{\SL}{\mathrm{SL}}
\DeclareMathOperator{\Sp}{\mathrm{Sp}}
\DeclareMathOperator{\Sym}{\mathbf{Sym}}

\renewcommand{\mod}{\mathop{\mathrm{mod}}\nolimits}

\renewcommand{\Im}{\mathop{\mathrm{Im}}\nolimits}

\def\ov{\overline}\newcommand{\car}[2]{\left[{\substack{#1\\#2}}\right]}

\def\carep{\boldsymbol{\varepsilon}}
\renewcommand{\theta}{\vartheta}

\def\diag{\textrm{diag}}

\newcommand{\M}{N}

\newcommand{\Thtri}[7]%
{\theta^{#7}\!\!\begin{bmatrix}#1#2#3\\#4#5#6\end{bmatrix}\!}

\begin{document}

\title{On a conjecture of Serre on abelian threefolds}

\author{Gilles Lachaud}
\address{Institut de Math\'ematiques de Luminy
\newline \indent
Universit\'e Aix-Marseille - CNRS
\newline \indent
Luminy Case 907, 13288 Marseille Cedex 9 - FRANCE
\newline \indent}
\email{lachaud@iml.univ-mrs.fr}

\author{Christophe Ritzenthaler}
\email{ritzent@iml.univ-mrs.fr}

\date{\today}

\keywords{Curve, Jacobian, abelian threefold, discriminant, Ciani quartic,
          modular form}

\subjclass[2000]{Primary 14G10; Secondary 14K25, 14H45}

\dedicatory{
\vskip5mm
\hfill \parbox{55mm}{En genre $3$, le th\'eor\`eme de Torelli s'applique
de fa\c{c}on moins satisfaisante : on doit extraire une myst\'erieuse
racine carr\'ee (J.-P.S., Collected Papers, n$^{o}$ 129)}}

\maketitle

\begin{abstract}
In this article, we give a reformulation of a result from Howe,
Leprevost and Poonen on a three dimensional family of abelian
threefolds. We also link their result to a conjecture of Serre on a
precise form of Torelli theorem for genus $3$ curves.
\end{abstract}

\tableofcontents


\newpage

\section{Introduction} \label{torelli}

\subsection{Geometric Torelli's theorem}
\label{torellig}

Let $K$ be an algebraically closed field. If $X$ is a (smooth
algebraic projective) curve of genus $g$ over $K$, the Jacobian
$\Jac X$ of $X$ is an abelian variety of dimension $g$, and $\Jac X$
has a canonical principal polarization $\lambda$. We obtain in this
way a morphism
$$
\begin{CD}
\Jac : \mathsf{M}_{g} & @>>> \mathsf{A}_{g}
\end{CD}
$$
from the space $\mathsf{M}_{g}$ of ($K$-isomorphism classes of)
curves of genus $g$ to the space $\mathsf{A}_{g}$ of
($K$-isomorphism classes of) $g$-dimensional principally polarized
abelian varieties (p.p.a.v.).

According to Torelli's Theorem, proved one century ago, the map $X \mapsto
(\Jac X, \Theta)$ is injective. An algebraic proof was provided by Weil
\cite{weil} half a century ago, and it is a long time studied question to
characterize the image of this map.

If $g = 3$, these spaces are both of dimension $3g-3=g (g+1)/2=6$.
According to Hoyt \cite{hoyt} and Oort and Ueno \cite{ueno}, the image of
$\mathsf{M}_{g}$ is exactly the space of indecomposable principally
polarized threefolds. Recall that $(A,\lambda)$ is decomposable if there is an
abelian subvariety $B$ of $A$ neither equal to $0$ nor to $A$, such
that the restriction of $\lambda$ to $B$ is a principal polarization, and
indecomposable otherwise. This was a problem left unsolved by Weil in
\cite{weil}.

Given a principally polarized abelian threefold $(A,\lambda)$ over
$K$, two natural questions arise :
\begin{enumerate}
\item How can we decide if the polarization is indecomposable ?
\item How can we decide if $A$ is the Jacobian of a a
hyperelliptic curve  ?
\end{enumerate}

Actually, both questions were answered by Igusa in 1967
\cite{igusa2} when $K = \CC$, making use of a particular modular
form $\chi_{18}$ on the Siegel upper half-space (see Th.
\ref{igusath} below).

\subsection{Arithmetic Torelli's theorem}

Assume now that $K$ is an arbitrary field. Then, as Serre noticed in
\cite{lauter}, the above correspondence is no longer true. Let
$(A,\lambda)$ be a p.p.a.v. of dimension $g$ over $K$, and assume
that $(A,\lambda)$ is isomorphic over $\ov{K}$ to the Jacobian of a
curve $\calX$ of genus $g$.

\begin{theorem}[Serre]
\label{twist}
The following alternative holds :
\begin{enumerate}
\item If $\calX$ is hyperelliptic, there exists a model $X/K$ of $\calX$ and a
$K$-isomorphism between the p.p.a.v. $(\Jac X,\Theta)$ and
$(A,\lambda)$.
\item If $\calX$ is non hyperelliptic, there exists a model
$X/K$ of $\calX$ and a quadratic character $\varepsilon :
\Gal(K_s/K) \to \{\pm 1\}$ such that $(\Jac X,\Theta)$ is
isomorphic to the twist $(A,\lambda)_{\varepsilon}$ of $(A,\lambda)$
by $\varepsilon$.
\end{enumerate}
In particular, if $\varepsilon$ is not trivial, this implies that
$\Jac X$ is not isomorphic to $A$ over $K$, but only over a
quadratic extension, and $(A,\lambda)$ is not isomorphic over $K$ to
the Jacobian of a curve. \qed

\end{theorem}

\subsection{Serre's conjecture}

Let us come back to the case $g = 3$. Let there be given an
indecomposable  principally polarized abelian threefold $(A,
\lambda)$ defined over $K$. In a letter to Top \cite{serre}  in
2003, J.-P. Serre asked two questions:
\begin{enumerate}
\item
How to decide, knowing only $(A, \lambda)$, that $X$ is hyperelliptic ?
\item
If $X$ is not hyperelliptic, how to find the quadratic character
$\varepsilon$ ?
\end{enumerate}
He proposed, in the case $K \subset \CC$, the following conjecture :

\begin{conjecture}
Let $(A,\lambda)$ be an undecomposable principally polarized abelian
threefold over $K$ isomorphic over $\ov{K}$ to the Jacobian of a
curve $\calX$ of genus $3$. Then there is an invariant
$\chi_{18}(A,\lambda)$ such that
\begin{enumerate}
\item
$\chi_{18}(A,\lambda) = 0$ is and only if $\calX$ is hyperelliptic;
\item
the character $\varepsilon$ is the one defined by the action of
$\Gal(\bar{K}/K)$ on the square root of $\chi_{18}(A,\lambda)$.
\end{enumerate}
\end{conjecture}
We use here the notation $\chi_{18}$ to emphasize that this
conjecture was inspired by the results obtained by Igusa, and also
much earlier by Klein \cite{klein} (see the remark after Cor.
\ref{kleincor}). In this article Klein relates (up to an undetermined
constant) the modular form $\chi_{18}$ and the square of the discriminant of
the quartic
$\calX$ (when $\chi_{18}(\tau) \ne 0$). This invariant seemed
to Serre a good choice to find this ``mysterious square root".

We plan to answer in the affirmative  this conjecture for a family
of abelian threefolds which are isogenous to the product of three
elliptic curves (see Cor. \ref{serrecor}). This will rely on the
work of Howe, Leprevost and Poonen \cite{leprevost} for which we
propose a natural rephrasing. For any field $K$ of characteristic
different from $2$, they consider abelian threefolds $(A,\lambda)$
defined as a quotient of three elliptic curves (with the trivial
polarization) by a certain subgroup of $2$-torsion points. For this
three-dimensional family, they make explicit the equation of the
related curve and express the character $\varepsilon$ by a invariant
$\mathsf{T}$ involving the coefficients of the elliptic curves. In
the first part, we show  that $\mathsf{T}$ can be naturally
interpreted as a determinant. In a second phase, we take $K \subset
\CC$ and by uniformization, we express $\mathsf{T}$ in terms of
certain Thetanullwerte of the elliptic curves. Then using the
duplication and transformation formula we express the modular form
$\chi_{18}(A,\lambda)$ in terms of the same Thetanullwerte and
compare the two expressions. We also obtain a proof of Klein's
result in this particular case and give the constant involved (see
Cor. \ref{kleincor}).

We describe now briefly the different sections. In Sec. \ref{CianiQuartics},
we define Ciani quartics, go back to the aforementioned results of
\cite{leprevost}, and show the relation with Serre's conjecture ($\S$
\ref{rel-serre}). In Sec. \ref{general}, we recall some general facts about
abelian varieties over $\CC$ (of arbitrary dimension) and introduce the
modular function $\chi_{k}$ ($\S$ \ref{chi-mod}). We prove Serre's conjecture
in Sec. \ref{comparison}. Finally, an appendix gathers some technical proofs,
in particular the modularity of the form $\chi_k$.\\

{\bf Acknowledgements.} We would like to thank J.-P. Serre for interesting discussions and S. Meagher for relevant remarks.

\section{Ciani Quartics}
\label{CianiQuartics}

In this section we reformulate a result of \cite{leprevost} on a
three-dimensional family of non-hyperelliptic genus $3$ curves. In
particular, this gives  a more natural point of view on Prop. 15 of
\cite{leprevost}.

\subsection{Definition of Ciani quartics}
Edgardo Ciani gave in 1899 \cite{ciani} a classification of nonsingular
complex plane quartics curves based on the number of involutions in their
automorphism group. We describe below the family of quartics admitting (at
least) two commuting involutions (different from identity).

Let $K$ be a field with $\Char K \neq 2$, and $\Sym_{3}(K)$ the
vector space of symmetric matrices of size $3$ with coefficients in
$K$. Let
$$
Q_{m}(x, y, z) = {^{t}\!v}.m.v,  \quad v = (x^2,y^2,z^2), \quad
m =
\begin{bmatrix}
a_{1} & b_{3} & b_{2} \\
b_{3} & a_{2} & b_{1} \\
b_{2} & b_{1} & a_{3}
\end{bmatrix}
\in \Sym_{3}(K).
$$
Then
$$
Q_{m}(x, y, z) = a_{1} x^{4} + a_{2} y^{4} + a_{3} z^{4} + 2 (b_{1}
y^{2} z^{2} + b_{2} x^{2} z^{2} + b_{3} x^{2} y^{2})
$$
is a ternary quartic, and the map $m \mapsto Q_{m}$ is an
isomorphism of $\Sym_{3}(K)$ to the vector space of ternary quartic
forms invariant under the three involutions
$$
\sigma_{1}(x,y,z) = (-x,y,z), \quad \sigma_{2}(x,y,z) = (x,-y,z), \quad
\sigma_{3}(x,y,z) = (x,y,-z).
$$
The form $Q_{m}$ is the zero locus of a plane quartic curve $X_{m}$,
whose automorphism group contains the Vierergruppe $V_{4} =
(\ZZ/2\ZZ)^{2}$.\\
If $X_{m}$ is a nonsingular curve, we  say that
$X_{m}$ is a \emph{Ciani quartic} and that $Q_{m}$ is a \emph{Ciani
form}. Now, E. Ciani (\textit{loc. cit.}) proved that a plane
quartic admitting two commuting involutions is geometrically
isomorphic to a Ciani quartic (a more recent reference is
\cite{STNB}).

\begin{proposition}
If $X$ is a  plane quartic curve defined over $K$, admitting at
least two commuting involutions, also defined over $K$, then there
is $m \in \Sym_{3}(K)$ such that $X$ is isomorphic to $X_{m}$ over
$K$.
\end{proposition}

\begin{proof}
Let $M_1$ and $M_2$ in $\PGL_3(K)$ inducing two commuting
involutions of $X$. Then $M_i^2 =\alpha_i \mathbf{I}$ with
$\alpha_{i} \in K$ and $\det(M_i)^2 = \alpha_i^3$, hence $\alpha_i$
is a square and we can assume, by dividing $M_i$ by
$\sqrt{\alpha_i}$, that $M_i^2 = \mathbf{I}$. The two matrices $M_i$
commute so we can diagonalize them in the same basis : after a
change of coordinates, we can suppose that $M_i$ are (projectively)
equal to
$$
I_1=
\begin{bmatrix} -1 & 0 & 0 \\ 0 & 1 & 0 \\ 0 & 0 & 1 \end{bmatrix}, \quad
I_2=
\begin{bmatrix} 1 & 0 & 0 \\ 0 & -1 & 0 \\ 0 & 0 & 1 \end{bmatrix}.
$$
This implies that a quartic equation $Q(x,y,z) = 0$ of $X$ in these new
coordinates must be invariant by the involutions $\sigma_{1}$ and $\sigma_{2}$
above, hence, $Q$ is a Ciani form.
\end{proof}

\subsection{Discriminant of a ternary form}

Our definition of a Ciani form includes that its zero locus must be
a nonsingular curve. This condition is fulfilled if and only if the
discriminant of the form is not $0$. In order to obtain a criterium
for this condition, we develop an algorithm for the discriminant of
a general ternary form.

The \emph{multivariate resultant} $\Res(f_{1},\dots,f_{n})$ of $n$
forms $f_{1}, \dots f_{n}$ in $n$ variables with coefficients in a
field $K$ is an irreducible polynomial in the coefficients of
$f_{1}, \dots f_{n}$ which vanishes whenever $f_{1}, \dots f_{n}$
have a common non-zero root. One requires that the resultant is
irreducible over $\ZZ$, \idest it has integral coefficients with
greatest divisor equal to $1$, and moreover
$$\Res(x_{1}^{d_{1}},\dots,x_{n}^{d_{n}}) = 1$$
for any $(d_1,\ldots,d_n) \in \NN^n$. The resultant exists and is
unique. There is a remarkable determinantal formula for the
resultant of $3$ ternary forms of the same degree $d$, due to
Sylvester; see \cite{GPZ} for a modern exposition and a proof. We
give this formula in the case $d = 3$. Then
$\Res(f_{1},f_{2},f_{3})$ is a form of degree $27$ in $30$ unknowns.
We shall express $\Res(f_{1},f_{2},f_{3})$ as the determinant of a
square matrix of size $15$.

Let $I_{d}$ be the set of sequences $\nu = (\nu_{1}, \nu_{2}, \nu_{3})$
with $\nu_{1} + \nu_{2} + \nu_{3} = d$. The monomials $x^{\nu} =
x_{1}^{\nu_{1}} x_{2}^{\nu_{2}} x_{3}^{\nu_{3}}$, where $\nu \in I_{d}$, form
the canonical basis of the space $V_{d}$ of ternary forms of degree
$d$, which is of dimension $(d + 1)(d + 2)/2$.

For any monomial $x^{\nu} \in V_{2}$ with $\nu_{1} + \nu_{2} + \nu_{3} = 2$,
we choose arbitrary representations
$$
f_{i} = x_{1}^{\nu_{1} + 1} f_{i,1} + x_{2}^{\nu_{2} + 1} f_{i,2} +
x_{3}^{\nu_{3}} f_{i,3} \quad (1 \leq i \leq 3),
$$
where $f_{i,j}$ are forms of degree $2 - \nu_{j}$, for $1
\leq i, j\leq 3$. Such a representation is always possible, although not
unique. Now we define
$$
S(x^{\nu}) = \det
\begin{bmatrix}
f_{1,1} & f_{1,2} & f_{1,3} \\
f_{2,1} & f_{2,2} & f_{2,3} \\
f_{3,1} & f_{3,2} & f_{3,3} \\
\end{bmatrix}
\, .
$$
Note that this determinant is indeed a ternary form of degree
$$(2 - \nu_{1}) + (2 - \nu_{2}) + (2 - \nu_{3}) = 6 - 2 = 4.$$
Since the monomials $x^{\nu}$ with $\nu \in I_{2}$ make up a basis of
$V_{2}$, we have thus defined a linear map
$S : V_{2}
\longrightarrow V_{4}$. We consider the linear map
$$
\begin{CD}
T : V_{1} \times V_{1} \times V_{1} \times V_{2} & @>>> & V_{4}
\end{CD}
$$
given by
$$
T(l_{1},l_{2},l_{3},g) = l_{1} f_{1} + l_{2} f_{2} + l_{3} f_{3} + S(g).
$$
Now one proves that the determinant of $T$ is independent of the choices
made in the definition of $S$, and \emph{Sylvester's formula} holds :
$$\Res(f_{1},f_{2},f_{3}) = \det T.$$
Generally speaking, the matrix of $T$ involves $864$ monomials. Now,
let $Q$ be a ternary form of degree $d$, and $X$ be the plane
projective curve which is the zero locus of $Q$. Call $q_{1}, q_{2},
q_{3}$ the partial derivatives of $Q$. The \emph{discriminant} of
$Q$ is $ \Disc Q = \Res(q_{1},q_{2},q_{3})$. It is a form of degree
$3(d - 1)^{2}$ in the coefficients of $Q$, and $X$ is non singular
if and only if $\Disc Q \neq 0$. The discriminant is an invariant of
ternary forms : if $g \in \GL_{3}(K)$, then
\begin{equation}
\label{repres}
\Disc (Q _\circ g) = (\det g)^{w} \, \Disc Q, \quad\quad \text{where} \
w  = d(d - 1)^{2}.
\end{equation}
If $Q$ is a quartic form  $\Disc Q$ is a form of degree $27$ in the coefficients, with $w = 36$ in
\eqref{repres}. Applying Sylvester's formula, we get :

\begin{proposition}
\label{DiscCiani} Let $m \in \Sym_{3}(K)$ and
$c_{i}=a_j a_k-b_i^2$ is the cofactor of $a_{i}$ for $1 \leq i \leq 3$. If $(q_{1},q_{2},q_{3})$
are the partial derivatives of the ternary quartic $Q_{m}$, then
$\Disc Q_m = 2^{54} \mathsf{D}(m)$, where
$$
\mathsf{D}(m) = a_{1}\,a_{2}\,a_{3} \, (c_{1}\,c_{2}\,c_{3})^{2} \,
\det(m)^4.\hfill \qed
$$

\end{proposition}

Note that this result was
obtained by Edge \cite{edge}, in a more intricate way.

We denote by $\mathsf{S}$ the set of $m \in \Sym_{3}(K)$ such that
$$
a_{1}a_{2}a_{3} \neq 0, \quad c_{1}c_{2}c_{3} \neq 0.
$$
Now, Prop.
\ref{DiscCiani} implies that the curve $X_m$ is nonsingular if and
only if $m$ belongs to the set $\mathsf{S}^{\times} = \mathsf{S}
\cap GL_{3}(K)$.

\begin{lemma}
\label{fromStoXm}
The map $m \mapsto Q_{m}$ from $\mathsf{S}^{\times}$ to the set
$\mathsf{Q}$ of Ciani forms is a bijection. \qed
\end{lemma}

The automorphisms of $X_m$ induce a  simple description of its
Jacobian. In order to make it explicit, we need to introduce a
certain product of elliptic curves.

\subsection{Product of elliptic curves}
\label{productec}

We introduce the following notations : let
$$
E_{i} : y^{2} = x(x^{2} - 4 b_{i} x - 4 c_{i}), \qquad (b_{i} \in K,
c_{i} \in K^{\times}) \quad (i = 1, 2, 3),
$$
three elliptic curves with $(0,0)$ as a rational $2$-torsion point.
The discriminant of $E_{i}$ is $\Delta_{i} = 2^{12} c_{i}^{2}
\delta_{i}$, where $\delta_{i} = b_{i}^{2} + c_{i} \in K^{\times}$.
We assume that there exists a square root $\rho \in K^{\times}$ of
$\delta(A) = \delta_{1}\delta_{2}\delta_{3}$, that is, $\Delta_{1}
\Delta_{2} \Delta_{3}$ is a square in $K$.  We denote by
$\mathsf{A}$ the set of products $E_1 \times E_2 \times E_3$ of such
curves and we define
$$\widetilde{\mathsf{A}} = \set{\widetilde{A} = (A,\rho) \in
\mathsf{A} \times K^{\times}} {\rho^{2} = \delta(A)}.$$ If
$\widetilde{A} \in \widetilde{\mathsf{A}}$, we put $a_{i} = \rho /
\delta_{i}$ and
$$
\mathbf{Mat}(\widetilde{A}) =
\begin{bmatrix}
a_{1} & b_{3} & b_{2} \\
b_{3} & a_{2} & b_{1} \\
b_{2} & b_{1} & a_{3}
\end{bmatrix} \, \in \mathsf{S}.
$$
Conversely, a matrix $m \in \mathsf{S}$ defines an abelian threefold
$A(m) \in \mathsf{A}$, which is the product of the curves
$$
E_{i} : y^{2} = x(x^{2} - 4 b_{i} x - 4 c_{i}) \qquad (i = 1, 2, 3).
$$
Then
$$
\delta_{i} = b_{i}^{2} + c_{i} = a_{j} a_{k} \in K^{\times}, \quad
\Delta_{1} \Delta_{2} \Delta_{3} =
(2^{18}\,a_{1}a_{2}a_{3}c_{1}c_{2}c_{3})^{2},\quad \delta(A) =
(a_{1}a_{2}a_{3})^{2}.$$ We define $\rho(m) = a_{1}a_{2}a_{3}$, and
$\mathbf{Ab}(m) = (A(m), \rho(m)) \in \widetilde{\mathsf{A}}$.

\begin{lemma}
\label{fromAtoS}
The maps
$$
\begin{CD}
\mathbf{Mat} & : \widetilde{\mathsf{A}} & @>>> & \mathsf{S}, &
\qquad & \mathbf{Ab}  & : \mathsf{S} & @>>> &
\widetilde{\mathsf{A}},
\end{CD}
$$
are mutually inverse bijections. \qed
\end{lemma}

The two lemmas \ref{fromStoXm} and \ref{fromAtoS} provide a natural
map from the set $\mathsf{Q}$ of Ciani quartic forms to
$\widetilde{\mathsf{A}}$. This map has actually a geometric meaning,
and in order to explain it, we introduce the following notation. If
$m \in \GL_{3}(K)$, we denote by $\Cof(m)$ the \emph{cofactor
matrix} of $m$, satisfying
$$
m.\tp{\Cof m}= \det m. \mathbf{I}, \quad \det \Cof m =
(\det m)^{2},\quad \Cof \Cof m = (\det m).m.
$$
Let $Q_m$ be a Ciani form associated to $m \in \mathsf{S}^{\times}$ and $X_m$  be the corresponding Ciani quartic
$$
X_{m} : Q_{m}(x, y, z) = F_{m}(x^{2}, y^{2}, z^{2}) = 0.
$$
By quotient, we get three genus one curves \begin{eqnarray*}
C_{1}:=  X_{m} / <1, \sigma_{1}> \; \;: F(y z, x^{2}, y^{2}) = 0, \\
C_{2} := X_{m} / <1, \sigma_{2}> \; \;: F(z x, y^{2}, z^{2}) = 0, \\
C_{3}:= X_{m} / <1, \sigma_{3}>  \; \; : F(x y, z^{2}, x^{2}) = 0,
\end{eqnarray*}
 where $\sigma_{i}$ $(i = 1,2,3)$
are the involutions of $X_m$. Another change of variables maps the
genus $1$ quartics $C_{i}$ to the elliptic curves
$$
F_{i} : y^{2} = x(x^{2} - 4 d_{i} x - 4 a_{i}\det(m)),
\qquad (i = 1,2,3).
$$
In this way, we get a map
$$
\varphi : X_{m} \longrightarrow B_{m} = F_{1} \times F_{2} \times F_{3}.
$$
Let us now look more closely at $B_{m}$. The identity
$$\Cof \Cof m = (\det m).m$$
implies that the cofactor of $c_{i}$ is $a_{i} \det m$. Hence,
$$
\mathbf{Ab}(\Cof m) = (B_{m}, c_{1}\,c_{2}\,c_{3}).
$$
Since the Jacobian is the Albanese variety of $X_{m}$, we get a
factorization
\medskip
$$
\begin{CD}
X_m & \\
@V{\iota}VV \quad \searrow ^{\varphi} \\
\Jac X_{m} & @>{\Phi}>> & B_{m}
\end{CD}
$$
\medskip

where $\iota$ is a canonical embedding. Since the images of the
regular differential forms on $F_{i}$ make a basis of those on $X_m$,
we obtain:

\begin{proposition}
The map
$$
\begin{CD}
\Phi : \Jac X_{m} & @>>> & A(\Cof m)
\end{CD}
$$
is a $(2,2,2)$-isogeny defined over $K$. \qed
\end{proposition}

The correspondences
$$
\begin{CD}
m         & @>>> & \Cof m \\
@VVV      &      & @VVV & \\
Q_{m}     & @>>> & A(\Cof m)  & @>{isg}>> & & \Jac X_{m}
\end{CD}
$$
lead to a commutative diagram, where $\mathsf{Q}$ is the space of
Ciani quartics over $K$:
$$
\begin{CD}
\mathsf{S}^{\times} & @>{\Cof}>>    & \mathsf{S}^{\times} \\
@VV{=}V             &               & @VV{\mathbf{Ab}}V & \\
\mathsf{Q}          & @>{``\Jac"}>> & \widetilde{\mathsf{A}}
\end{CD}
$$
In the next section we describe the kernel of the isogeny $\Phi$.

\subsection{The theory of Howe, Leprevost and Poonen revisited}
\label{revisited}

The previous isogeny can be made more precise as we recall from
\cite{leprevost}. There are some differences between their notation
and ours, see the remark at the end of Sec.\ref{expli-diff} for a comparison.
Let us introduce couples $(A, W)$, where:
\begin{enumerate}
\item
$A \in \mathsf{A}$ as defined in $\S$ \ref{productec}.
\end{enumerate}
The Weil pairings on the factors combine to give a non degenerate
alternating pairing $e_{2}$ on the finite group scheme $A[2]$ over
$K$.
\begin{enumerate}
\setcounter{enumi}{1}
\item
$W$ is a totally isotropic indecomposable subspace of $A[2]$ defined
over $K$.
\end{enumerate}
Choose a basis $(P_{i},Q_{i}) \in E_{i}[2]$, that is, a level $2$
structure on $E_{i}$. This defines a level $2$ structure on $A$.
In \cite[Lem.13]{leprevost} it is proved that after a labeling of
the $2$-torsion points we can write
\begin{equation}
\label{DefW}
W=
\begin{Bmatrix}
(O,O,O), & (O, Q_{2}, Q_{3}), & (Q_{1}, O, Q_{3}), & (Q_{1}, Q_{2}, O), \\
(P_{1},P_{2},P_{3}), & (P_{1},R_{2},R_{3}), & (R_{1},P_{2},R_{3}), &
(R_{1},R_{2},P_{3})
\end{Bmatrix}
\end{equation}
with
$$
Q_{i} = (0,0), \quad P_{i} = (0, 2 b_{i} + \rho_{i}\,), \quad R_{i}
= (0, 2 b_{i} - \rho_{i}\,), \quad \rho_{1} \rho_{2} \rho_{3} =
\rho_{W},
$$
and the four possible choices of $\rho_{1}, \rho_{2}, \rho_{3}$ leading to
the same value of $\rho_{1} \rho_{2} \rho_{3}$ give the same subgroup $W$.
Conversely, if $(A, \rho) \in \widetilde{\mathsf{A}}$ is given, if we
choose $\rho_{1}, \rho_{2}, \rho_{3}$ in such a way that $\rho_{1} \rho_{2}
\rho_{3} = \rho$, and if we define $P_{i}$ and $R_{i}$ as above, then we can
define a subgroup $W_{\rho}$ by \eqref{DefW}.

\begin{lemma}
The map $(A, \rho) \mapsto (A, W_{\rho})$ from
$\widetilde{\mathsf{A}}$ to the set of couples $(A, W)$ as defined
above is a bijection. \qed
\end{lemma}

We take on $A = A(m)$ the principal polarization $\lambda$ which is the
product of the canonical polarizations on each factor. Then we have
a commutative diagram
$$
\begin{CD}
A     & @>{2 \lambda}>> & A^{\vee} \\
@V{\pi}V{d^{\circ}8}V  & & @A{d^{\circ}8}A{\widehat{\pi}}A \\
A'     & @>{\lambda'}>> & (A')^{\vee} \\
\end{CD}
$$
with a unique principal polarization $\lambda'$ on $A' = A'(m) =
A(m)/W_{\rho(m)}$. From \mbox{\cite[Prop.15]{leprevost}} we get :

\begin{theorem} \label{HLPrevisited}
\label{Jacobienne} The composition of the isogeny $\Phi$ and of the
projection $\pi$ leads to an isomorphism of p.p.a.v.:
$$
\begin{CD}
\Jac X_{m} & @>>> & A'(\Cof m).
\end{CD}
\rlap \qed
$$
\end{theorem}

As a corollary, for any $m$, the p.p.a.v. $(A', \lambda')$ is
indecomposable.

The reverse direction is more interesting and will give an algebraic
answer to Serre's conjecture.

\subsection{Relation with Serre's conjecture}
\label{rel-serre}

We need the following elementary lemma from linear algebra.

\begin{lemma}
\label{cofactor}
The map $m \mapsto \Cof m$ induces an exact sequence
$$
\begin{CD}
1 & @>>> \{\pm 1\} & @>>> \GL_{3}(K) & @>{\Cof}>> G^{\times 2}(K) &
@>>> 1
\end{CD}
$$
with
$$
G^{\times 2}(K) = \set{m\in \GL_{3}(K)}{\det m \in K^{\times 2}}.
\rlap \qed
$$
\end{lemma}

Let $(A,\rho)=\widetilde{A} \in \widetilde{\mathsf{A}}$ and
$m=\mathbf{Mat}(\widetilde{A})$. Denote
$$
\mathsf{T}(\widetilde{A}) := \det(m) = 2 b_{1} b_{2} b_{3} -
\rho(\frac{b_{1}^2}{\delta_{1}} + \frac{b_{2}^2}{\delta_{2}} +
\frac{b_{3}^2}{\delta_{3}} - 1).
$$

\begin{theorem}
\label{howebis}
The following results hold.
\begin{enumerate}
\item if $\mathsf{T}(\widetilde{A}) = 0$, that is, $m \in \mathsf{S}
\setminus \mathsf{S}^{\times}$, there is a
hyperelliptic curve $X$ of genus $3$ such that $A'(m)$ is isomorphic
to the Jacobian of $X$.
\item if $\mathsf{T}(\widetilde{A}) \ne 0$, that is, $m \in
\mathsf{S}^{\times}$, then there exists a non hyperelliptic curve of
genus $3$ defined over $K$ whose Jacobian is isomorphic to $A'(m)$
if and only if $\mathsf{T}(\widetilde{A})$ is a square in $K$.
\end{enumerate}
\end{theorem}

\begin{proof}
The first part is \cite[Prop.14]{leprevost} where the
hyperelliptic curve is constructed explicitly.
For the second part, if $\det(m)$ is a square, then using Lem.
\ref{cofactor}, we see that there exists a matrix $m' \in
\mathsf{S}^{\times}$ such that $m=\Cof(m')$ and we apply Th.
\ref{Jacobienne}.
If $d=\det(m)$ is not a square, let $m_d=d m$ and $\widetilde{A}_{d}
= \mathbf{Ab}(m_{d})$.  $A(m_d)$ is defined by
$$
E_{i} : y^{2} = x(x^{2} - 4 b_{i} d x - 4 c_{i} d^{2}) \qquad (i =
1, 2, 3),
$$
Thus $A(m_d)$ is a quadratic twist of $A(m)$. Now, $\det(m_d)$ is
a square, so there exists $m'$ such that $\Jac(X_{m'})$ is
isomorphic to $A'(m_d)$. Since $A'(m_d)$ is a quadratic twist of
$A'(m)$ and is the Jacobian of a non hyperelliptic curve, Th.
\ref{twist} shows that $A'(m)$ cannot be a Jacobian.
\end{proof}

\begin{corollary}
With the same notation as above:
\begin{enumerate}
\item
if $\mathsf{T}(\widetilde{A}) \in K^{\times 2}$, there is an isogeny
defined over $K$
$$
\begin{CD}
\Jac X_{m'} & @>>> & A(m), \quad \Cof m' = m,
\end{CD}
$$
\item
If $\mathsf{T}(\widetilde{A}) \notin K^{\times 2}$, there is an
isogeny defined over $K$
$$
\begin{CD}
\Jac X_{m'} & @>>> & A(m_{d})  \quad \Cof m' = d m.
\end{CD}
\rlap \qed
$$
\end{enumerate}
\end{corollary}

We hope to give in a near future a geometric interpretation of the
connection of Serre's problem with the determinant of certain quadratic
forms in the general case.

\begin{remark*}
\label{expli-diff}
In \cite{leprevost}, Howe, Leprevost and Poonen  write the elliptic
curves
$$
y^{2} = x(x^{2} + A_{i} x + B_{i}) \qquad \text{avec} \ A_{i}, B_{i}
\in K \qquad (i = 1, 2, 3).
$$
So
$$
A_{i} = - 4 b_{i}, \quad B_{i} = - 4 c_{i},
$$
$$
\Delta_{i} = A_{i}^{2} - 4 B_{i} = 16 (b_{i}^{2} + c_{i}) = 16
\delta_{i},
$$
$$
d_{i} = -(A_{i} + 2 x(P_{i})) = 4 b_{i} - 4(b_{i} + \rho_{i}) = - 4
\rho_{i}, \quad d_{i}^{2} = \Delta_{i}.
$$
And the factor
$$
T_{0}(\widetilde{A}) = d_{1} d_{2} d_{3}
(\frac{A_{1}^{2}}{\Delta_{1}} + \frac{A_{2}^{2}}{\Delta_{2}} +
\frac{A_{3}^{2}}{\Delta_{3}}  - 1) - 2 A_{1} A_{2} A_{3},
$$
which is
$$
T_{0}(\widetilde{A}) =  64 [- \rho (\frac{b_{1}^{2}}{\delta_{1}} +
\frac{b_{2}^{2}}{\delta_{2}} + \frac{b_{3}^{2}}{\delta_{3}} - 1) + 2
b_{1}b_{2}b_{3}].
$$
We then have $T_{0}(\widetilde{A}) = 64 \, \mathsf{T}(\widetilde{A})$.
\end{remark*}

\section{Complex abelian varieties}
\label{general}

We recall in this section some well known propositions on abelian
varieties over $\CC$ and fix the notation.

\subsection{The symplectic group}
If $V$ is a module of rank $2g$ over a commutative ring $R$ and if $E$ is a
non-degenerate  alternating bilinear form on $V$, a basis $(a_{i})_{1 \leq i
\leq 2g}$ of $V$ is said \emph{symplectic} if the matrix $(E(a_{i},a_{j})) =
J$, where
$$
J =
\begin{bmatrix} 0 & \mathbf{1}_{g} \\ - \mathbf{1}_{g} & 0 \end{bmatrix} \, .
$$
The group of matrices
$$
M = \begin{bmatrix}
A & B \\ C & D \end{bmatrix} \in \SL_{2g}(R)
$$
such that  $M.J.\tp{M} = J$ is the \emph{symplectic group} $\Sp_{2g}(R)$. It
acts simply transitively on the set of symplectic bases of $V$.
\begin{lemma}[{\cite[Lem.8.2.1]{lange}}] \label{eqcondi}
If $M \in \Sp_{2g}(R)$ the following conditions
are equivalent.
\begin{enumerate}
\item
$M \in \Sp_{2g}(R)$.
\item
$\tp{A}.C$ and $\tp{B}.D$ are symmetric, and $\tp{A}.D - \tp{C}. B =
\mathbf{1}_g$.
\item
$A \tp B$ and $C \tp D$ are symmetric and $A .\tp{D} - B.\tp{C} =
\mathbf{1}_g$. \qed
\end{enumerate}
\end{lemma}
The group $\Sp_{2g}(R)$ is the group of $R$-rational points of a
Chevalley group scheme $\Sp_{2g}$, which contains certain remarkable subgroups
defined as follows. The reductive subgroup $\mathbf{M}$ of $\Sp_{2g}$ is the
subgroup which respects the canonical decomposition $\ZZ^{2g} = \ZZ^{g}
\oplus \ZZ^{g}$.  Elements of $\mathbf{M}(\ZZ)$ are
$$
M =
\begin{bmatrix}
A & 0 \\ 0 & \tp{A}^{-1}
\end{bmatrix} \, ,
\quad  A \in \GL_{g}(\ZZ).
$$
The unipotent subgroup $\mathbf{U}$ is the stability group leaving
pointwise fixed the canonical totally isotropic subspace $V_{0}$, which is
the first direct summand in the standard decomposition. Elements of
$\mathbf{U}(\ZZ)$ are
$$
U =
\begin{bmatrix}
\Id_{g} & B \\ 0 & \Id_{g}
\end{bmatrix} \, , \quad \tp{B} = B,  \quad B \in \mathbf{M}_{g}(\ZZ).
$$
The unipotent subgroup $\mathbf{V}$ opposite to $\mathbf{U}$ is the
stability group of the second direct summand in the standard
decomposition. One has $\mathbf{V} = \tp{\,\mathbf{U}} =
J.\mathbf{U}.J^{-1}$.  Elements of $\mathbf{V}(\ZZ)$ are
$$
V =
\begin{bmatrix}
\Id_{g} & 0 \\ C & \Id_{g}
\end{bmatrix} \, , \quad \tp{C} = C,  \quad C \in \mathbf{M}_{g}(\ZZ).
$$
The subgroup $\mathbf{P} = \mathbf{M} \ltimes \mathbf{U}$ is the
parabolic subgroup of $\Sp_{2g}$ normalizing $V_{0}$, and $\mathbf{P}$ is
actually a maximal parabolic subgroup. Elements of $\mathbf{P}(\ZZ)$ are
$$
P =
\begin{bmatrix}
A & B \\ 0 & \tp{A}^{-1}
\end{bmatrix} \, , \quad A.\tp{B} = B.\tp{A}, \quad  A \in \GL_{g}(\ZZ),
\quad B \in \mathbf{M}_{g}(\ZZ).
$$

\subsection{Abelian varieties}

Let $\Omega = [w_{1} \dots w_{2g}] \in \mathbf{M}_{g,2g}(\CC)$,
where $w_{1}, \dots, w_{2g}$ are columns vectors giving a basis of
$\CC^g$ on $\RR$. It generates a lattice
$$\Lambda = \Omega \ZZ^{2g} \subset \CC^{g}.$$
Let  $\EuScript{R}_{g}$ be the set of matrices $\Omega \in
\mathbf{M}_{g,2g}(\CC)$ satisfying  the \emph{Riemann conditions}
$$
\Omega.J. \tp{\Omega} = 0, \quad
2i (\overline{\Omega}.J^{-1}.\tp{\Omega})^{-1} > 0
$$
($> 0$ means positive definite). We call such a  matrix $\Omega$ a
\emph{period matrix}. If $\Omega \in \EuScript{R}_{g}$, the torus $A_{\Omega} =
\CC^{g} / \Lambda$ is an abelian variety of dimension $g$ with a principal
polarization $\lambda$ represented by the hermitian form $H =
2i (\overline{\Omega}.J^{-1}.\tp{\Omega})^{-1}$ (see \cite[Lem.4.2.3]{lange}).\\
The group $\GL_{g}(\CC)$ acts on the left on $\EuScript{R}_{g}$. If we
write
$$
\Omega = [(w_{1} \, \dots \, w_{g}) \, (w_{g + 1}  \, \dots \,
w_{2g})] = \Rh{\Omega_{1}}{\Omega_{2}}, \quad \textrm{where} \
\Omega_{i} \in \mathbf{M}_{g}(\CC),
$$
we get $W.\Rh{\Omega_{1}}{\Omega_{2}} =
\Rh{W.\Omega_{1}}{W.\Omega_{2}}$ for any $W \in \GL_{g}(\CC)$.  This
action induces  an isomorphism of p.p.a.v. In
particular if we choose $W = \Omega_2^{-1}$, we see that $A_{\Omega}$ is
isomorphic to the p.p.a.v.
$$
A_{\tau} = A_{\boldsymbol{\Omega}(\tau)}, \quad
\boldsymbol{\Omega}(\tau) = \Rh{\tau}{\mathbf{1}_{g}}, \quad  \tau =
\boldsymbol{\tau}(\Omega) = \Omega_{2}^{- 1}\Omega_{1},
$$
and $\Omega \in \EuScript{R}_{g}$ if and only if  $\boldsymbol{\tau}(\Omega)$
belongs to the Siegel upper half plane
$$
\HH_{g} = \set{\tau \in \mathbf{M}_{g}(\CC)}{\tp{\tau} = \tau, \ \Im
\tau > 0}.
$$
We call a matrix $\tau \in \HH_{g}$ a \emph{Riemann matrix}. The \emph{Siegel
modular group}
 $\Gamma_{g} = \Sp_{2g}(\ZZ)$ acts
on the right on  $\EuScript{R}_{g}$: if $\Omega \in \EuScript{R}_{g}$ and
if $M \in \Gamma_{g}$,
$$
\Omega.M = \Rh{\Omega_{1}}{\Omega_{2}} \begin{bmatrix} A & B \\
C & D \end{bmatrix} = \Rh{\Omega_{1}A + \Omega_{2}C}{\Omega_{1}B +
\Omega_{2}D}.
$$
This action corresponds to a change of symplectic basis. The group
$\Gamma_{g}$ also acts on the left on the Siegel upper half plane :
if $\tau \in \HH_{g}$, we denote $$M. \tau = (A \tau + B)(C \tau +
D)^{-1}.
$$
Both actions are linked by
\begin{equation}
\label{TauTransp} M.\boldsymbol{\tau}(\Omega) =
\boldsymbol{\tau}(\Omega. \tp M).
\end{equation}

\subsection{Isotropy and quotients}

For any maximal isotropic subgroup $V \subset \FF_{2}^{2g}$, we have the
\emph{transporter}
$$
\Trans(V) = \set{M \in \Sp_{2g}(\FF_{2})}{M V_{0} = V},
$$
$V_0$ being the canonical maximal isotropic subgroup generated by the
vectors $e_1,\ldots,e_g$ of the canonical basis. Since $\Sp_{2g}(\FF_2)$
permutes transitively the maximal isotropic subgroups of $\FF_{2}^{2g}$, the
transporter is a left coset: $\Trans(V) = M_0  \mathbf{P}(\FF_2)$,
for any $M_0 \in \Trans(V)$. Hence, the set of maximal isotropic
subgroups is the quotient set $\Sp_{2g}(\FF_2)/\mathbf{P}(\FF_2)$, a set with
$135$ elements if
$g = 3$.

Let now $\Omega \in \EuScript{R}_{g}$, $\Lambda=\Omega \ZZ^{2g}$ and
$(A,\lambda) = (A_{\Omega},H)$ be the corresponding p.p.a.v. of dimension $g$.
The linear map $\alpha : \ZZ^{2g}
\rightarrow \tfrac{1}{2} \Lambda$ such that
$$
\alpha(x) = \frac{1}{2}\Omega.x
$$
defines a level $2$ symplectic structure on $A[2]$, that is, an
isomorphism
$$
\begin{CD}
\bar{\alpha} : \FF_{2}^{2g} & @>{\sim}>> A[2]
\end{CD}
$$
and if $V \subset \FF_{2}^{2g}$ is a maximal isotropic subgroup, the same
property holds for $W = \bar{\alpha}(V) \subset A[2]$. If $\pi :\CC^{3}
\rightarrow
\CC^{3}/ \Lambda$ is the canonical projection, the lattice $\Lambda_{W} =
\pi^{-1}(W)$ is associated to $A/W$ as the following diagram shows:
$$
\begin{CD}
0 & @>>> & \Lambda_{W} / \Lambda & @>>> &
\CC^{3} / \Lambda & @>>> & \CC^{3}/ \Lambda_{W} & @>>> 0 \\
& & & & @| & & @| & & @| \\
0 & @>>> & W & @>>> & A  & @>>> & A/ W & @>>> 0.
\end{CD}
$$
 We define
$$\Trans(W) = \set{M \in \Gamma_g}{M (\mod 2) \in \Trans(V)}.$$
We introduce now the congruence subgroup
\begin{eqnarray*}
\Gamma_{0,g}(2) & = & \set{
\begin{bmatrix}
A & B \\ C & D
\end{bmatrix}
\in \Gamma_{g}}{C \equiv 0 (\mod 2)}.
\end{eqnarray*}
This is the transposed subgroup of the group $\Gamma_{g}^{0}(2)$, see
$\S$ \ref{modularity}. From Prop. \ref{Generators} we deduce
that $\Gamma_{0,g}(2) = \mathbf{P}(\ZZ).\Gamma_g(2)$, hence
$$\Trans(W) = M \Gamma_{0,g}(2)$$
for any $M \in  \Trans(W)$.

\begin{proposition}
\label{actiondemi}
With the previous notation, if $\tau = \boldsymbol{\tau}(\Omega)$ then
$\frac{1}{2} \tp M. \tau$ is a Riemann matrix of the p.p.a.v. $A/W$ for all $M
\in \Trans(W)$.
\end{proposition}

\begin{proof}
If $M \in \Trans(W)$, $W \mod{\Lambda}$ is generated by the vectors
$$ \frac{1}{2} \, \Omega.M e_{1}, \quad \dots \quad \frac{1}{2} \,
\Omega.M e_{g},
$$
and the matrix $$\Omega'= \Omega.M.H, \qquad H =
\begin{bmatrix}
\frac{1}{2} \mathbf{1}_{g} & 0 \\ 0 & \mathbf{1}_{g}
\end{bmatrix} \, ,
$$
generates $\Lambda_{W}$. Using \eqref{TauTransp}, we  get
$$
\boldsymbol{\tau}(\Omega')=\tp(M.H).\tau=\frac{1}{2} \tp M \tau.
$$
By \cite[Prop. 16.8]{milne}, the polarization $2\lambda$ of $A$
reduces to a principal polarization $\lambda'$ on $A'=A/W$. This last corresponds
canonically to $\Omega'$ since
$$
2i (\ov{\Omega'}. J^{-1}.  \tp{\Omega}')^{-1}=
2i (\ov{\Omega }  M H J^{-1} \tp H \tp M \tp \Omega)^{-1} =
2  \cdot 2i (\ov{\Omega}  J^{-1}  \tp \Omega)^{-1}.
$$
\end{proof}

\subsection{Theta functions}

We recall the definition of theta functions
with (entire) characteristics $[\carep] = \car{\varepsilon_1}{\varepsilon_2}$ where
$\varepsilon_1,\varepsilon_2 \in \ZZ^g$, following \cite{lange}. The
\emph{(classical) theta function} is
$$
\theta  \car{\varepsilon_1}{\varepsilon_2}(z,\tau)=
\sum_{n \in \mathbb{Z}^g} q^{(n+\varepsilon/2) \tau (n+\varepsilon/2) + 2
(n+\varepsilon/2)(z+\varepsilon_2/2)} \quad
(\tau \in \HH_g, z \in \CC^g).
$$
The \emph{Thetanullwerte} are the values at $z = 0$ of these functions, and
we denote
$$
\theta \car{\varepsilon_1}{\varepsilon_2}(\tau) =
\theta \car{\varepsilon_1}{\varepsilon_2}(0,\tau).
$$



We now state two formulas.

\begin{proposition}[duplication formula, see {\cite[Cor.IIA2.1]{rauch} and
\cite[IV.th.2]{igusa1}}]
Let $\car{\varepsilon_1}{\varepsilon_2}$ and $\car{\varepsilon_1}{\delta}$
be two characteristics and $\tau \in \HH_g$.  Then
\begin{equation} \label{dupliform}
\theta \car{\varepsilon_1}{\varepsilon_2}(\tau/2) \theta
\car{\varepsilon_1}{\delta}(\tau/2)  =
\sum_{\mu \in(\mathbb{Z}/2 \mathbb{Z})^g}
(-1)^{\mu \delta} \cdot \theta
\car{ \varepsilon_1-\mu}{\varepsilon_2-\delta}(\tau)  \cdot \theta
\car{\mu}{\varepsilon_2-\delta}( \tau).
\end{equation}

\end{proposition}

The second is called \emph{transformation formula}.
 Let $M=\left(\begin{array}{cc} A & B
\\ C & D \end{array}\right) \in  \Sp_{2g}(\mathbb{Z})$. We let $M$
acts on the characteristics in the following way
$$[M.\carep]=M.\car{\varepsilon_1}{\varepsilon_2}=\car{D \varepsilon_1 - C \varepsilon_2+ (C \tp D)_0}{
-B \varepsilon_1+ A \varepsilon_2 + (A \tp B)_0}$$ where $P_0$ denotes the
diagonal of the matrix $P$.

\begin{proposition}[{\cite[V.\S.2]{igusa1}}] \label{transfo}
\begin{equation*}
\theta [M.\carep](M. \tau)= \kappa(M) \cdot w^{
\phi_{[\varepsilon_1,\varepsilon_2]}(M)} \cdot j(M,\tau)^{1/2} \cdot
\theta [\carep](\tau)
\end{equation*}
where $\kappa(M)^2$is a root of $1$ depending only on $M$, $\omega=e^{i \pi/4}$,
$$j(M,\tau)=\det(C \tau+D)$$
 and
$$\phi_{[\varepsilon_1,\varepsilon_2]}(M)=\varepsilon_1 \tp D B \varepsilon_1-2 \varepsilon_1 \tp B C \varepsilon_2+
\varepsilon_2 \tp C A \varepsilon_2-  2 ( D \varepsilon_1- C \varepsilon_2)
\cdot (A \tp B)_0. \rlap \qed$$

\end{proposition}

We will need a slightly modified version of the
previous result.

\begin{corollary} \label{ctransfo1}
For any characteristic $\car{\varepsilon_1'}{\varepsilon_2'}$ and for any
$M \in \Sp_{2g}(\ZZ)$ we have
\begin{equation}  \label{formulatransfo}
\theta  \car{\varepsilon_1'}{\varepsilon_2'} (M. \tau)=
c(M,\tau) \cdot \omega^{-\phi_{[\varepsilon_1,\varepsilon_2]}(M^{-1})} \cdot \theta
\car{\tp A (\varepsilon_1'- (C \tp D)_0)+ \tp C
(\varepsilon_2'- (A \tp B)_0) }{\tp B (\varepsilon_1'- (C \tp D)_0) + \tp
D (\varepsilon_2'- (A \tp B)_0)} (\tau)
\end{equation}
where $$c(M,\tau)=\kappa(M^{-1})^{-1} \cdot j(M,\tau).$$
\end{corollary}

\begin{proof}
To inverse the action on the characteristics, we let
$\tau'=M^{-1}. \tau$ in the transformation formula. Note that
$$
M^{-1}=\left(\begin{array}{cc} \tp D & -\tp B
\\ -\tp C & \tp A \end{array}\right)
$$
and that
$[M. \carep]=[\tp M^{-1} \binom{\varepsilon_1}{
\varepsilon_2}]+ \car{(C \tp D)_0}{(A \tp B)_0}.$
Thus we get the action on the characteristics. For the factor $j(M,\tau)$ note that
\begin{eqnarray*}
j(M,\tau)&=&\det(C \tau+D)=\det(C M^{-1}. \tau'+D) \\
&=& \det(C (\tp D \tau'-\tp B)
(-\tp C \tau'+\tp A)^{-1}+D) \\
&=& \det(C (\tp D \tau'-\tp B)+ D(-\tp C \tau'+\tp A)) \det(-\tp C \tau'+\tp A)^{-1}\\
&=& \det(-\tp C \tau'+\tp A)^{-1}=j(M^{-1},\tau')^{-1}
\end{eqnarray*}
using Lem. \ref{eqcondi}.
\end{proof}

\begin{corollary} \label{ctransfo2}
Let $\Omega=\Rh{\Omega_1}{\Omega_2}$ be a period matrix and
$\tau=\boldsymbol{\tau}(\Omega)=\Omega_2^{-1} \Omega_1 \in \HH_g$.
Let $\Omega'=\Omega \tp M=\Rh{\Omega_1'}{\Omega_2'}$. Then
$$j(M,\boldsymbol{\tau}(\Omega))=\det(\Omega_2)^{-1} \cdot  \det(\Omega_2').$$
\end{corollary}
\begin{proof}
We compute
\begin{eqnarray*}
\det(C \tau+D)& =&  \det(\tau \tp C+ \tp D) \\
&=& \det(\Omega_2)^{-1} \det(\Omega_1 \tp C+ \Omega_2 \tp D) \\
&=&  \det(\Omega_2)^{-1} \cdot  \det(\Omega_2'),
\end{eqnarray*}
the last expression coming from (\ref{TauTransp}).
\end{proof}

\subsection{The modular function $\chi_k$}
\label{chi-mod}

Recall that a characteristic $\car{\varepsilon_1}{\varepsilon_2}$ is
\emph{even} (resp. \emph{odd}) if $\varepsilon_1.\varepsilon_2
\equiv 0 \pmod{2}$ (resp. $\varepsilon_1.\varepsilon_2 \equiv 1
\pmod{2}$). Let $S_g$ (resp. $U_g$) be the set of even (resp. odd)
characteristics $\car{\varepsilon_1}{\varepsilon_2}$ with coefficients in
$\{0,1\}$. It is well known that
$$\# S_g = 2^{g-1}(2^g+1), \; \; \#U_g=2^{g-1}(2^g-1).$$
Let $\Omega= \Rh{\Omega_1}{\Omega_2} \in \EuScript{R}_{g}$ and
$\tau=\Omega_2^{-1} \Omega_1 \in \HH_g$ be a Riemann matrix. For $g
\geq 2$, we denote $k=\# S_g/2$ and we are interested in the
following expressions :
$$\chi_k(\tau)=\prod_{\carep \in S_g} \theta [\carep](\tau).
$$


Recall that  a function $f$ is a \emph{modular form} of weight $w$ for the congruence
subgroup $\Gamma \in \Gamma_g$ if for all $\tau \in \HH_g$ and $M \in \Gamma$ one has
$$f(M.\tau)=j(M,\tau)^w f(\tau).$$

Using Cor.\ref{ctransfo2}, we get
\begin{corollary}
Let $f$ be a modular form of weight $k$  for $\Gamma$ on $\HH_g$. For $\Omega=
\Rh{\Omega_1}{\Omega_2} \in \EuScript{R}_{g}$, we define
$\tau=\Omega_2^{-1} \Omega_1 \in \HH_g$  a Riemann matrix
and $$f(\Omega):=\det(\Omega_2)^{-k} \cdot f(\tau).$$
Then for all $M \in \Gamma$
$$f(\Omega. M)=f(\Omega). \rlap \qed$$
\end{corollary}

In his beautiful paper \cite{igusa2}, Igusa proves the following result
[\textit{loc. cit.}, Lem. 10 \& 11]. Denote by $\Sigma_{140}$ the
thirty-fifth elementary symmetric function of the eighth power of the
even Thetanullwerte.

\begin{theorem}
\label{igusath}
For $g \geq 3$, the product $\chi_k(\tau)$ is a modular form of weight $k$ for the
group $\Gamma_g$. Moreover, If $g = 3$ and $\tau \in \HH_3$, then:
\begin{enumerate}
\item
$A_{\tau}$ is decomposable if $\chi_{18}(\tau) = \Sigma_{140}(\tau) = 0$.
\item
$A_{\tau}$ is a hyperelliptic Jacobian if $\chi_{18}(\tau) = 0$ and
$\Sigma_{140}(\tau) \neq 0$.
\item
$A_{\tau}$ is a non hyperelliptic Jacobian if $\chi_{18}(\tau) \neq 0$.
\qed
\end{enumerate}
\end{theorem}

This theorem gives an answer to the two questions raised in Sec.\ref{torellig} over
$\CC$.\\

In the sequel, we will need the following result to prove the independence of our
results from the choices we will make. The proof is the case $g = 3$ of Th.
\ref{Modularity}.

\begin{proposition}
\label{modulardemi}
The product $\tau \mapsto \chi_{18}(\tfrac{1}{2} \tau)$ is a modular form on $\HH_{3}$ of
weight $18$ for $\Gamma_{3}^{0}(2)$.
\end{proposition}

\section{Comparison of analytic and algebraic discriminants}
\label{comparison}

In this part, we make the link between the algebraic result Th.\ref{howebis}
and Serre's conjecture on the modular function $\chi_{18}$. To do so, we first
compute a quantity related easily to $T(\widetilde{A})$ in terms of the
Thetanullwerte on the elliptic curves. Then, after a good choice of a
symplectic matrix $N$ (related to the subgroup $W$ we use for the quotient),
we compute $\chi_{18}((\tp N. \tau)/2)$ in terms of the same
Thetanullwerte. Thus, we express $\chi_{18}$ on the quotient $A/W$. Finally
 we compare the expressions to get Serre's conjecture.


\subsection{Expression of the algebraic discriminant}

We come back to the hypotheses of $\S$ \ref{productec}, and
specialize to the case $K \subset \CC$. Let $A = E_{1} \times E_{2}
\times E_{3} \in \mathsf{A}$, where
$$
E_{i} : y^{2} = x(x^{2} - 4 b_{i} x - 4 c_{i}), \qquad (b_{i} \in K,
c_{i} \in K^{\times}) \quad (i = 1, 2, 3).
$$
We choose a root $\rho$ of $\delta(A)$, and put $m = \mathbf{Mat}(\widetilde{A})$
 with $\widetilde{A}=(A,\rho)$. Let
\begin{equation*}
\mathsf{X}(m) : =  (a_{1}\,a_{2}\,a_{3})^{4} \,
(c_{1}\,c_{2}\,c_{3})^{2} \, \det m.
\end{equation*}
Since $\mathsf{T}(\widetilde{A})=\det m$,
$\mathsf{T}(\widetilde{A})$ is a square in $K$ if and only if $X(m)$
is a square in $K$. The function $\mathsf{X}$ appears naturally in
our problem since it is related to the function $\mathsf{D}(m)$
(Prop. \ref{DiscCiani}) by
\begin{equation} \label{XD}
\mathsf{X}(\Cof m) = \mathsf{D}(m)^{2},
\end{equation}
and this reflects Serre's conjecture according to Th.\ref{HLPrevisited}.
In order to determine the expression of $X(m)$ is terms of the Thetanullwerte,
we use the following uniformization.
The curves $E_{i}$ can be written as
$$
E(\omega_{1i},\omega_{2i}) : y^{2} = x (x + \frac{\pi^2}{\omega_2^2}
\theta \car{0}{0}(\tau_i)^4) (X+\frac{\pi^2}{\omega_2^2} \theta
\car{0}{1}(\tau_i)^4),
$$
with $\Rh{\omega_{1i}}{\omega_{2i}} \in \EuScript{R}_{1}$ and $\tau_i=\frac{\omega_{1i}}{\omega_{2i}}$. We identify
$\EuScript{R}_{1}^{3}$ with the set of matrices
$$
\Omega = \Rh{\Omega_{1}}{\Omega_{2}} = \left[ \
\begin{pmatrix}
\omega_{11} & 0           & 0           \\
0           & \omega_{12} & 0           \\
0           & 0           & \omega_{13}
\end{pmatrix}
\begin{pmatrix}
\omega_{21} & 0           & 0          \\
0           & \omega_{22} & 0          \\
0           & 0           & \omega_{23}
\end{pmatrix}
\right]
$$
such that
$$\tau = \mathbf{\tau}(\Omega)=\Omega_{2}^{- 1}\Omega_{1} =
\begin{bmatrix}
\tau_{1} & 0 & 0 \\
0 & \tau_{2} & 0 \\
0 & 0 & \tau_{3}
\end{bmatrix}
\in \HH_3 .
$$
We define
$$
A(\Omega) := E(\omega_{11},\omega_{21}) \times
E(\omega_{12},\omega_{22}) \times E(\omega_{13},\omega_{23}),
$$
$$
\rho(\Omega): =
\frac{\pi^{6}}{64(\det \Omega_{2})^{2}} \,
\Thtri{1}{1}{1}{0}{0}{0}{4}(\tau).
$$
This defines an element $\widetilde{A}(\Omega): = (A(\Omega),\rho(\Omega)) \in
\widetilde{\mathsf{A}}$, and a matrix $m(\Omega): =
\mathbf{Mat}(\widetilde{A}(\Omega))$. For $1 \leq i \leq 3$, denote
$$
\theta_{0i} = \theta \car{0}{0}(\tau_i), \quad
\theta_{1i} = \theta \car{1}{0}(\tau_i), \quad
\theta_{2i} = \theta \car{0}{1}(\tau_i).
$$
The coefficients of $A(\Omega)$ and $m(\Omega)$ are
\begin{eqnarray*}
a_{i} & = &
- \frac{\pi^{2}}{4} \,
\frac{\omega_{2i}^{2}}{(\omega_{2j} \omega_{2k})^{2}} \,
\frac{\theta_{1j}^{4} \theta_{1k}^{4}} {\theta_{1i}^{4}} \, ,\\
b_{i} & = & -
\frac{\pi^{2}}{4 \omega_{2i}^{2}} (\theta_{0i}^{4} + \theta_{2i}^{4}), \\
c_{i} & = & - \frac{\pi^{4}}{4 \omega_{2i}^{4}} \,
\theta_{0i}^{4} \, \theta_{2i}^{4},
\end{eqnarray*}
where $(i,j,k)$ is a cyclic permutation.
The determinant of $m(\Omega)$ is expressed as follows. Let
$$
\mathsf{a} = \theta_{01}^2 \theta_{02}^2 \theta_{23}^2, \quad
\mathsf{b} = \theta_{01}^2 \theta_{22}^2 \theta_{03}^2, \quad
\mathsf{c} = \theta_{21}^2 \theta_{02}^2 \theta_{03}^2, \quad
\mathsf{d} = \theta_{21}^2 \theta_{22}^2 \theta_{23}^2,
$$
$$
\mathsf{R}_{1} =
(\mathsf{a} + \mathsf{b} + \mathsf{c} + \mathsf{d})
(\mathsf{a} + \mathsf{b} - \mathsf{c} - \mathsf{d})
(\mathsf{a} - \mathsf{b} - \mathsf{c} + \mathsf{d})
(\mathsf{a} - \mathsf{b} + \mathsf{c} - \mathsf{d}),
$$
Then
$$
\det m(\Omega) = \frac{\pi^6}
{2^4 \cdot \prod_{i=1}^3 (\omega_{2i}^2 \cdot
(\theta_{0i}^4-\theta_{2i}^4))}
\cdot \mathsf{R}_{1}.
$$
Thus we get
\begin{eqnarray*}
\mathsf{X}(m(\Omega))
& = &
(\frac{\pi^{12}}{2^6 \cdot \prod_{i=1}^3
\omega_{2i}^4} \cdot \prod_{i=1}^3 \theta_{0i}^4 \theta_{2i}^4)^2 \cdot
(\frac{\pi^6}{2^6 \cdot \prod_{i=1}^3 \omega_{2i}^2} \cdot
\prod_{i=1}^3 (\theta_{0i}^4-\theta_{2i}^4))^4 \\
&  & \cdot
( \frac{\pi^6}{2^4  \cdot   \prod_{i=1}^3 (\omega_{2i}^2 \cdot
(\theta_{0i}^4-\theta_{2i}^4))} \cdot \mathsf{R}_{1}) \\
& = & \frac{\pi^{54}}{2^{40}} \cdot \det(\Omega_2)^{-18} \cdot(
\prod_{i=1}^3 \theta_{0i}^8 \theta_{2i}^8
(\theta_{0i}^4-\theta_{2i}^4)^3) \cdot \mathsf{R}_{1}.
\end{eqnarray*}

\subsection{The subgroup $W$}
\label{DefN}
With the notation of Sec. \ref{revisited}, we can always assume the following correspondences
$$
\quad P_{i} \leftrightarrow
\frac{\omega_{1i}}{2}, \quad Q_{i} \leftrightarrow
\frac{\omega_{2i}}{2}, \quad R_{i} = P_{i} + Q_{i}
\leftrightarrow\frac{\omega_{1i} + \omega_{2i}}{2}
$$
for the points of $E_i$. The characteristics associated to the
points of $W$ (see $\S$ \ref{revisited}) are
$$
\begin{array}{cccccccc}
\car{0 0 0}{0 0 0}, & \car{0 0 0}{0 1 1}, & \car{0 0 0}{1 0 1}, &
\car{0 0 0}{1 1 0}, & \car{1 1 1}{ 0 0 0}, & \car{1 1 1}{0 1 1}, &
\car{1 1 1}{1 0 1}, & \car{1 1 1}{1 1 0}.
\end{array}
$$
It defines a maximal isotropic subgroup $V$ of $\FF_2^6$. A basis of
$V$ over $\FF_{2}$ is given by the three vectors
$$
\begin{array}{cccccccc}
\alpha_{1} = \car{0 0 0}{0 1 1}, & \alpha_{2} = \car{0 0 0}{ 1 1 0},
& \alpha_{3} = \car{1 1 1}{0 0 0}
\end{array} \, .
$$
The matrix
$$
\M =
\begin{bmatrix}
0   &   0 &   1 &   0 & - 1 & 0\\
0   &   0 &   1 &   0 &   0 & 0\\
0   &   0 &   1 & - 1 &   0 & 0\\
0   &   1 &   0 &   0 &   0 & 0 \\
- 1 & - 1 &   0 &   0 &   0 & 1 \\
  1 &   0 &   0 &   0 &   0 & 0
\end{bmatrix}
$$
belongs to $\Gamma_{3}$ and satisfies $\M.e_{i} \equiv \alpha_{i} \pmod{2} \
\text{if} \ 1 \leq i \leq 3$, thus $\M \in \Trans(W)$.\\
The set
\begin{eqnarray*}
\Gamma_{g}(1,2) & = & \set{
\begin{bmatrix}
A & B \\ C & D
\end{bmatrix}
\in \Gamma_{g}}
{(A.\tp{B})_{0} \equiv (C.\tp{D})_{0} \equiv 0 (\mod 2)},
\end{eqnarray*}
 is a subgroup of $\Gamma_{g}$, and $\kappa^{2}$ is a
character of $\Gamma_{g}(1,2)$ \cite[p. 181]{igusa1}.

\begin{lemma}
\label{kappa}
The matrices $\M$ and $\tp \M$ are in $\Gamma_{3}(1,2)$, and
$\kappa(\M)^{2} = \kappa(\tp \M)^{2} = \pm 1$.
\end{lemma}

\begin{proof}
We have $\M = LQ$, where
$$
L = \begin{bmatrix}
A & 0 \\ 0 & \tp{A}^{-1}
\end{bmatrix}
, \quad \text{with} \ A =
\begin{bmatrix}
0  & -1 & 1 \\
0  &  0 & 1 \\
-1 & 0  & 1
\end{bmatrix},
$$
and
$$
Q =
\begin{bmatrix}
0   &   0 &   0 &   1 &   0 & 0 \\
0   &   0 &   0 &   0 &   1 & 0 \\
0   &   0 &   1 &   0 &   0 & 0 \\
- 1 &   0 &   0 &   0 &   0 & 0 \\
0   & - 1 &   0 &   0 &   0 & 0 \\
0   &   0 &   0 &   0 &   0 & 1
\end{bmatrix} \, .
$$
\bigskip

One checks easily that $L, \tp{L}, Q, \tp{Q}$ belong to $\Gamma_{3}(1,2)$,
hence, $\M$ and $\tp{\M}$ are in $\Gamma_{3}(1,2)$ as well. If
$$
M =
\begin{bmatrix}
A & B \\ 0 & D
\end{bmatrix}
\in \mathbf{P}(\ZZ),
$$
then $\kappa(M)^{2} = \det D$, see \cite[Lem. 7, p. 181]{igusa2}. Now
$$
Q^{2} =
\begin{bmatrix}
S & 0 \\ 0 & S
\end{bmatrix}
, \quad \text{with} \ S =
\begin{bmatrix}
-1 & 0  & 0 \\
0  & -1 & 0 \\
0  & 0  & 1
\end{bmatrix}.
$$
From this we deduce that
$$
\kappa(L)^{2} = \det A=1, \quad \kappa(Q)^{4} = \kappa(Q^{2})^{2} = \det S=1,
$$
hence, $\kappa(\M)^{2} = \kappa(\tp{\M})^{2} = \pm 1$.
\end{proof}

\begin{proposition}
Let $\Omega'= \Omega \M H$. Then
$$\mathbf{\tau}(\Omega')=\frac{1}{2} \tp \M. \tau$$
is
a Riemann matrix for $A'(m)$.
Moreover, the value $\chi_{18}(\Omega')$ is independent on the choice of $N
\in \Trans(W)$.
\end{proposition}

\begin{proof}
The first assertion comes from Prop. \ref{actiondemi}, the second from
Prop. \ref{modulardemi}.
\end{proof}

\subsection{Expression of $\chi_{18}(\Omega')$ as a discriminant}

Our main result in this section is the following

\begin{theorem}
\label{KleinIdentity}
Let $\Omega \in \EuScript{R}_{1}^{3}$ and $A(\Omega)$ be the corresponding
abelian threefold, let $m = m(\Omega)\in \mathsf{S}$ be the associated
matrix, and $\Omega'\in
\EuScript{R}_{3}$ be a Riemann matrix of $A(\Omega)/W$. Then
$$
\left(\frac{\pi}{2}\right)^{54} \cdot \chi_{18}(\Omega')
= \mathsf{X}(m).
$$
\end{theorem}

\begin{proof}
The strategy is the following. Let $N$ be the matrix defined in $\S$
\ref{DefN}, and define $\tau' = \tp{\M}.\tau = 2 \mathbf{\tau}(\Omega')$.
\begin{enumerate}
\item
\label{strategy1}
Pair the  Thetanullwerte in $\tau'/2$ such that one
can apply the duplication formula \eqref{dupliform}. We then obtain expressions
in terms of  Thetanullwerte in $\tau'$. Such a pairing is not unique
and one makes here a choice which allows an easy comparison of the final
formulas.
\item
\label{strategy2}
For each of the  Thetanullwerte in $\tau'$, apply the transformation
formula \eqref{formulatransfo} to obtain an expression in $\tau$.
\item
\label{strategy3}
Finally, since $\tau=\diag(\tau_1,\tau_2,\tau_3)$, we get
$$
\theta \car{a_1b_1c_1}{a_2b_2c_2} (\tau)=\prod_{i=1}^3
\theta \car{a_1}{b_1}(\tau_i).
$$
\end{enumerate}

Let
$$c(\M)=\kappa(\tp \M^{-1})^{-2} \det(\Omega_2)^{-1} \det(\Omega_2')
=
\pm \det(\Omega_2)^{-1} \det(\Omega_2')$$
by  Lem.
\ref{kappa}.\\
 Applying steps \eqref{strategy1} to
\eqref{strategy3} with the software MAGMA (see
\url{http://iml.univ-mrs.fr/~ritzenth/programme/check2.m}), we get
the following 18 identities, where we write
$$
\theta \car{0 0 0}{0 0 0} \theta \car{0 0 0}{0 0 1} =
\theta \car{0 0 0}{0 0 0}(\tau'/2) \theta \car{0 0 0}{0 0
1}(\tau'/2),  \quad c = c(\M).
$$
We make the pairing in such a way that the expressions of $\theta \car{0 0
0}{\varepsilon_2} \theta \car{0 0 0}{\delta}$ do not
contain $\theta_{1i}$ terms. The first four are, with the preceding notation,
\begin{eqnarray*}
\theta \car{0 0 0}{0 0 0} \theta \car{0 0 0}{0 0 1}
& = & c \, (\mathsf{a} + \mathsf{b} + \mathsf{c} + \mathsf{d}) \\
\theta \car{0 0 0}{0 1 0} \theta \car{0 0 0}{0 1 1}
& = & c \, (\mathsf{a} + \mathsf{b} - \mathsf{c} - \mathsf{d})\\
\theta \car{0 0 0}{1 0 0} \theta \car{0 0 0}{1 0 1}
& = & -c \, (\mathsf{a} - \mathsf{b} - \mathsf{c} + \mathsf{d}) \\
\theta \car{0 0 0}{1 1 0} \theta \car{0 0 0}{1 1 1} & = & -c \,
(\mathsf{a} - \mathsf{b} + \mathsf{c} - \mathsf{d})
\end{eqnarray*}
and the remaining $14$ are
\begin{eqnarray*}
\theta \car{0 1 0}{0 0 0} \theta \car{0 1 0}{0 0 1}
& = &
2 c \,
(\theta_{01}\theta_{21}\theta_{02}\theta_{22}\theta_{03}^2 +
 \theta_{01}\theta_{21}\theta_{02}\theta_{22}\theta_{23}^2) \\
\theta \car{1 0 0}{0 0 0} \theta \car{1 0 0}{0 0 1}
& = & 2 c \,
(\theta_{01}^2 \theta_{02}\theta_{22}\theta_{03}\theta_{23} +
 \theta_{21}^2 \theta_{02}\theta_{22}\theta_{03}\theta_{23}) \\
\theta \car{1 1 0}{0 0 0} \theta \car{1 1 0}{0 0 1}
& = & 2c \,
(\theta_{01}^2\theta_{21}\theta_{02}^2\theta_{03}\theta_{23} +
\theta_{01}\theta_{21}\theta_{22}^2\theta_{03}\theta_{23})\\
\theta \car{0 1 0}{1 0 0} \theta \car{0 1 0}{1 0 1}
& = & 2 c \,
(\theta_{01}^2\theta_{21}\theta_{02}\theta_{22}\theta_{03}^2  -
 \theta_{01}^2\theta_{21}\theta_{02}\theta_{22}\theta_{23}^2) \\
\theta \car{1 0 0}{0 1 0} \theta \car{1 0 0}{0 1 1}
& = & 2  c \,
(\theta_{01}^2\theta_{02}\theta_{22}\theta_{03}\theta_{23}  -
 \theta_{21}^2\theta_{02}\theta_{22}\theta_{03}\theta_{23}) \\
\theta \car{1 1 0}{1 1 0} \theta \car{1 1 0}{1 1 1}
& = & -2 c \,
(\theta_{01}\theta_{21}\theta_{02}^2\theta_{03}\theta_{23} -
 \theta_{01}\theta_{21}\theta_{22}^2\theta_{03}\theta_{23}) \\
\theta \car{0 0 1}{0 0 0} \theta \car{0 0 1}{0 1 0}
& = & 2 c \,
(\theta_{01}\theta_{11}\theta_{02}\theta_{12}\theta_{03}\theta_{13}) \\
\theta \car{0 0 1}{1 0 0} \theta \car{0 0 1}{1 1 0}
& = & 2c \,
(\theta_{01}\theta_{11}\theta_{02}\theta_{12}\theta_{03}\theta_{13}) \\
\theta \car{0 1 1}{1 1 0} \theta \car{0 1 1}{1 0 0}
& = & 2 c \,
(\theta_{11}\theta_{21}\theta_{12}\theta_{22}\theta_{03}\theta_{13}) \\
\theta \car{1 0 1}{0 0 0} \theta \car{1 0 1}{0 1 0}
& = & 2 c \,
(\theta_{01}\theta_{11}\theta_{12}\theta_{22}\theta_{13}\theta_{23}) \\
\theta \car{1 1 1}{0 0 0} \theta \car{1 1 1}{1 1 0}
& = & 2  c \,
(\theta_{11}\theta_{21}\theta_{02}\theta_{12}\theta_{13}\theta_{23}) \\
\theta \car{0 1 1}{0 1 1} \theta \car{0 1 1}{1 1 1}
& = & -2 c \,
(\theta_{11}\theta_{21}\theta_{12}\theta_{22}\theta_{03}\theta_{13}) \\
\theta \car{1 1 1}{0 1 1} \theta \car{1 1 1}{1 0 1}
& = & -2 c \,
(\theta_{11}\theta_{21}\theta_{02}\theta_{12}\theta_{13}\theta_{23}) \\
\theta \car{1 0 1}{1 0 1} \theta \car{1 0 1}{1 1 1}
& = & -2 c \,
(\theta_{01}\theta_{11}\theta_{12}\theta_{22}\theta_{13}\theta_{23})
\end{eqnarray*}

Denote by $\mathsf{R}_{1}'$ the product of the first four lines. Obviously
$\mathsf{R}_{1}'=c(\M)^4 \mathsf{R}_{1}$. Calling $\mathsf{R}_{2}'$ the
product of the last fourteen lines, we get
$$\mathsf{R}_{2}'= 2^{14}  \cdot c(\M)^{14} \cdot  \left( \prod_{i=1}^3
\theta_{0i}^8 \theta_{2i}^8 (\theta_{0i}^4-\theta_{2i}^4)^3\right).$$
So
\begin{multline*}
\chi_{18}(\tau'/2) =
 \mathsf{R}_{1}' \mathsf{R}_{2}'= 2^{14} \cdot c(\M)^{18} \cdot
\left(\frac{2^{40}}{\pi^{54}} \cdot \det(\Omega_2)^{18} \right)
\cdot \mathsf{X}(m) \\ =
\left(\frac{2}{\pi}\right)^{54} \cdot  \det(\Omega_2')^{18} \cdot  \mathsf{X}(m),
\end{multline*}
which is the expected result.
\end{proof}

Since $\mathsf{X}(m)$ is equal to $\mathsf{T}(\widetilde{A})$ up to a square, Th.\ref{howebis}
and Th.\ref{KleinIdentity} show Serre's conjecture.

\begin{corollary} \label{serrecor}
Let $K \subset \CC$ and $m \in \mathsf{S}^{\times}$ with
coefficients in $K$. Let $A'(m)$ be the associated abelian threefold
and $\Omega'$ be one of its period matrix. Then
$$
\left(\frac{\pi}{2}\right)^{54} \cdot \chi_{18}(\Omega') \in
K^{\times 2}
$$
if and only if $A'(m)$ is the Jacobian of a non hyperelliptic genus $3$ curve.
\qed
\end{corollary}

In other words, Serre's conjecture is true for our three dimensional family
$\mathsf{A}$ of abelian threefolds.
\begin{corollary}
\label{kleincor}
If $m \in \mathsf{S}^{\times}$ and $\Omega_m$ is a period matrix
associated to the non hyperelliptic genus $3$ curve $X_m$ with Ciani form $Q_m$ then
$$
\chi_{18}(\Omega_{m})=\left(\frac{1}{2 \pi}\right)^{54} \cdot \Disc (Q_m)^2. \rlap \qed$$
\end{corollary}
\begin{proof}
Using  Th.\ref{HLPrevisited}  and \eqref{XD} we get
\begin{eqnarray*}
\left(\frac{\pi}{2}\right)^{54} \cdot \chi_{18}(\Omega_m) &=& \mathsf{X}(\Cof m) \\
&=& \mathsf{D}(m)^2 = (2^{-54} \cdot \Disc Q_m)^2.
\end{eqnarray*}
\end{proof}

\begin{remark*}
When $m \in \mathsf{S} \setminus \mathsf{S}^{\times}$, the abelian variety $A'(m)$
comes from a hyperelliptic curve and the above formula degenerates. However
in \cite{lockhart} and \cite{guardia} we find a beautiful
formula for the hyperelliptic case in every genus. Let
$$
C : Y_2=a_{2g+2} X^{2g+2}
+ \ldots + a_0= a_{2g+2} (X-\alpha_1) \cdots (X-\alpha_{2g+2})
$$
and
$$
\Delta_{\textrm{alg}}(C)=a_{2g+2}^{4g+2} \prod_{j<k}
(\alpha_j-\alpha_k)^2.
$$
They define also a modular form on $\HH_g$
$$\delta(\tau)=\prod_{\carep \in T} \theta[\carep](\tau)^8$$
where $T$ is a certain subset of even theta characteristic. One has
$$\Delta_{\textrm{alg}}(C)^{2n}=(2\pi)^{4rg}
\det(\Omega_1)^{-4r} \delta(\tau)^2$$
where
$$r = \binom{2g+2}{g+1}, \quad n = \binom{ 2g}{g+1},$$
and $\tau=\mathbf{\tau}(\Omega)$ for a certain period matrix
$\Omega=[\Omega_1,\Omega_2]$ of $\Jac(C)$.
\end{remark*}

\begin{remark*}
\label{kleinformula}
Denote by $V_{3}^{4}$ the $15$-dimensional affine open set of ternary
quartics.  Felix Klein proved in 1889 that there is a map
$$
\Omega : V_{3}^{4} \longrightarrow \EuScript{R}_{g}
$$
such that if $\Omega(Q) = \Rh{\Omega_{1}}{\Omega_{2}}$ and $X : Q=0$, then
$\Jac X = A_{\Omega(Q)}$ and
$$
\chi_{18}(\Omega) = c \Disc(Q)^{2},
$$
with some unspecified constant $c \in \CC$. We prove here that $c =
(1/2 \pi)^{54}$. Using this precise version of Klein's formula, it
is almost obvious to extend our theorem to the general case.
However, we did not include it, for we think that a good
presentation should include a modern proof of Klein's result. We
plan to do this in a forthcoming article.
\end{remark*}

\appendix
\label{appendix}
\section{}

\subsection{Modularity of $\chi_{k}$}
\label{modularity}
Let
$$
\Gamma_g(2) = \set{M \in \Gamma_g}{M \equiv \mathbf{1}_{2g} (\mod 2)}
$$
an recall that the sequence
$$1 \to \Gamma_g(2) \to \Gamma_g \to \Sp_{2g}(\FF_2) \to 1$$
is exact. We introduce the congruence subgroup
\begin{eqnarray*}
\Gamma_{g}^{0}(2) & = & \set{
\begin{bmatrix}
A & B \\ C & D
\end{bmatrix}
\in \Sp_{g}(\ZZ)}{B \equiv 0 (\mod 2)}.
\end{eqnarray*}
We need a set of generators for this subgroup. For any integer $n
\geq 1$, define
$$
M(n) = \mathbf{M}(\ZZ) \cap \Gamma_{g}(n), \quad
U(n) = \mathbf{U}(\ZZ) \cap \Gamma_{g}(n), \quad
V(n) = \mathbf{V}(\ZZ) \cap \Gamma_{g}(n).
$$
with
$$
\Gamma_g(n) = \set{M \in \Gamma_g}{M \equiv \mathbf{1}_{2g} (\mod n)}.
$$

\begin{proposition}
\label{Generators}
The subgroups $M(1)$, $U(2)$ and $V(1)$ generate $\Gamma_{g}^{0}(2)$, and
$\Gamma_{g}^{0}(2) = \Gamma_{g}(2).\mathbf{M}(\ZZ).\mathbf{V}(\ZZ)$.
\end{proposition}

\begin{proof}
First, the subgroups $M(2)$, $U(2)$ and $V(2)$ generate $\Gamma_{g}(2)$,
see \cite[p. 179]{igusa1}. Let
\begin{eqnarray*}
\Gamma_{g}^{1}(2) & = & \set{
\begin{bmatrix}
A & B \\ C & D
\end{bmatrix}
\in \Sp_{g}(\ZZ)}
{A \equiv D \equiv 1 (\mod 2) \ \text{and} \ B \equiv 0 (\mod 2)}.
\end{eqnarray*}
There is the following diagram, where the vertical arrow is the transpose
of the reduction modulo $2$:
$$
\begin{array}{ccccccc}
\Gamma_{g}(2) & \subset & \Gamma_{g}^{1}(2) & \subset & \Gamma_{g}^{0}(2) &
\subset & \Gamma_{g}(1) \\
\downarrow && \downarrow && \downarrow &&\downarrow \\
0 & \subset & \mathbf{U}(\FF_{2}) & \subset &
\mathbf{P}(\FF_{2}) & \subset & \Sp_{g}(\FF_{2}) \\
\end{array}
$$

Then, if $M \in \Gamma_{g}^{1}(2)$ is written as
usual
$$
\begin{bmatrix}
A & B \\ C & D
\end{bmatrix}
\begin{bmatrix}
\Id_{g} & 0 \\ C & \Id_{g}
\end{bmatrix}
=
\begin{bmatrix}
A + B C & B \\ C + D C & D
\end{bmatrix}
\in \Gamma_{g}(2),
$$
and if $M \in \Gamma_{g}^{0}(2)$, then
$$
\begin{bmatrix}
A & B \\ C & D
\end{bmatrix}
\begin{bmatrix}
A^{-1} & 0 \\ 0 & \tp{A}
\end{bmatrix}
=
\begin{bmatrix}
\Id_{g} & B\tp{A} \\ C A^{-1} & D\tp{A}
\end{bmatrix}
\in \Gamma_{g}^{1}(2).
$$
\end{proof}

\begin{theorem}
\label{Modularity}
Assume $g \geq 3$.
The function
$\chi_{k}(\tfrac{1}{2} \tau)$ is a modular form on $\HH^{g}$ of weight $k$
for $\Gamma_{g}^{0}(2)$.
\end{theorem}

\begin{proof}
J.-I. Igusa proved \cite[p. 850]{igusa2} that $\chi_{k}(\tau)$ is a
modular form of weight $k$ for $\Gamma_{g}(1)$ if $g \geq 3$. Let
$$
H =
\begin{bmatrix}
\frac{1}{2} \Id_{g} & 0 \\ 0 & \Id_{g}
\end{bmatrix} \, , \quad H.\tau = \tfrac{1}{2} \tau.
$$
Let $f(\tau) = \chi_{k}(\tfrac{1}{2} \tau) = \chi_{k}(H.\tau)$. It is
sufficient to check that
$$
f(M.\tau) = j(M,\tau)^{k} f(\tau)
$$
if $M$ belongs to one of the generating subgroups described in Prop.
\ref{Generators}. First, if $M \in M(1)$, then $H.M = M.H$, hence,
$$
f(M.\tau) = \chi_{k}(H.M.\tau) = \chi_{k}(M.H.\tau) =
j(M,H.\tau)^{k} \chi_{k}(H.\tau) = j(M,\tau)^{k} f(\tau),
$$
since $j(M,\tau) = \pm 1$ for every $M \in \mathbf{M}(\ZZ)$ does not depend on $\tau
\in \HH_{g}$. Now, if $U \in U(2)$, then
$$
U = U'^{2} =
\begin{bmatrix}
\Id_{3} & 2 B \\ 0 & \Id_{3}
\end{bmatrix}
\, , \quad \text{where} \ U' =
\begin{bmatrix}
\Id_{3} & B \\ 0 & \Id_{3}
\end{bmatrix}
\in \mathbf{U}(\ZZ),
$$
and $H.U = H.U'^{2} = U'.H$. This implies
$$
f(U.\tau) = \chi_{k}(H.U'^{2}.\tau) = \chi_{k}(U'.H.\tau) =
j(U',H.\tau)^{k} \chi_{k}(H.\tau) = j(U,\tau)^{k} f(\tau),
$$
since $j(U^{n},\tau) = 1$ for every $U \in \mathbf{U}(\ZZ)$. If $V \in
V(1)$, then $H.V = V^{2}.H$. Hence
\begin{multline*}
f(V.\tau) = \chi_{k}(H.V.\tau) = \chi_{k}(V^{2}.H.\tau) =
j(V^{2},H.\tau)^{k} \, \chi_{k}(H.\tau) \\= j(V,\tau)^{k} \chi_{k}(H.\tau)
= j(V,\tau)^{k} f(\tau),
\end{multline*}
since $j(V^{2},\tau) = j(V,2\tau)$ for every $V \in V(1)$.
\end{proof}



\begin{thebibliography}{99}

\bibitem{STNB}
Bars, Francesc, Automorphism groups of genus $3$ curves, in \textit{Corbes de
g\`enere 3}, Notes del Seminari de Teoria de Nombres de Barcelona {\bf 14}, 2006,
27-62.

\bibitem{lange}
C. Birkenhake, , H. Lange,  \textit{Complex abelian varieties} Second
edition. Grundleh\-ren der Mathematischen Wissenschaften, {\bf 302}
Springer-Verlag, Berlin, 2004.

\bibitem{ciani}
E. Ciani,  I Varii Tipi Possibili di Quartiche Piane pi\`u Volte
Omologico-Armoniche, Rend. Circ. Mat. Palermo {\bf 13} (1899), 347-373.

\bibitem{edge}
W.L. Edge, The discriminant of a certain ternary quartic, Proc. Roy.
Soc. Edinburgh, Sect. A. {\bf 62} (1948), 268-272.

\bibitem{guardia}
J. Gu\`ardia,  Jacobian nullwerte and algebraic equations. J. Algebra {\bf 253}
(2002), 112-132.

\bibitem{GPZ}
I.M. Gel'fand, M.M. Kapranov, A.V. Zelevinsky
\textit{Discriminants, resultants, and multidimensional determinants},
Birkh\"{a}user, Boston, (1994).

\bibitem{hoyt} W.L. Hoyt,  On products and algebraic families of Jacobian varieties.  Ann. of Math. {\bf 77},  (1963), 415-423.

\bibitem{leprevost}
E. Howe, F. Lepr\'evost, B. Poonen,  Large torsion subgroups of split
Jacobians of curves of genus two or three. Forum Math. {\bf 12}, (2000),  315-364.

\bibitem{igusa2}
J.-I. Igusa,  Modular forms and projective invariants,
Amer. J. Math, {\bf 89}, (1967), 817-855.

\bibitem{igusa1}
J.-I. Igusa, \textit{ Theta functions}, Grundlehren der
mathematischen Wissenschaften, {\bf 194}, Springer Verlag, (1972).

\bibitem{klein}
F. Klein,  Zur Theorie der Abelschen Funktionen. Math. Annalen,  {\bf 36}
(1889-90) = Gesammelte mathematische Abhandlungen,  {\bf  XCVII},
388-474.

\bibitem{lauter}
K. Lauter,  Geometric methods for improving the upper bounds on the
number of rational points on algebraic curves over finite fields. With an
appendix by Jean-Pierre Serre. J.~Algebraic Geom. {\bf 10}, (2001), 19-36.

\bibitem{lockhart}
P. Lockhart,  On the discriminant of a hyperelliptic curve. Trans. Amer. Math.
Soc.  {\bf 342}, (1994), 729-752.

\bibitem{milne}
J.S. Milne,  Abelian varieties, in \textit{ Arithmetic geometry} (Storrs, Conn., 1984), 103-150,
Springer, New York, (1986).

\bibitem{ueno} F. Oort, K. Ueno,  Principally polarized abelian varieties of dimension two or three are Jacobian varieties.  J. Fac. Sci. Univ. Tokyo Sect. IA Math., {\bf  20},  (1973), 377-381.

\bibitem{rauch}
H.E. Rauch, H.M. Farkas, \textit{ Theta functions with applications to Riemann
surfaces} The Williams \& Wilkins Co., Baltimore, Md., (1974).

\bibitem{serre}
J.-P. Serre, Letter to Jaap Top, February 28, private communication, (2003).


\bibitem{weil}
A. Weil,  Zum Beweis des Torellischen Satzes. Nachr. Akad. Wiss.
G\"{o}ttingen. Math.-Phys. Kl. IIa. (1957), 33-53 ; = \OE uvres Sc., vol.
II, [1957a], 307-327, Springer, New York, (1979).

\end{thebibliography}
\end{document}